\documentclass[journal]{IEEEtran}

\usepackage{graphicx}
\usepackage{amsmath,amssymb}

\usepackage{array}
\usepackage[noadjust]{cite}
\usepackage[justification=centering, font=scriptsize]{caption}

\begin{document}

\title{A Prosumer-Centric Framework for Concurrent Generation and Transmission Planning — Part I}

\author{Ni~Wang,
        Remco~Verzijlbergh,
        Petra~Heijnen,
        and~Paulien~Herder
\thanks{N. Wang, R. Verzijlbergh and P. Heijnen is with the Faculty of Technology, Policy and Management, Delft University of Technology, Delft, the Netherlands (e-mails: \{n.wang, r.a.verzijlbergh, p.w.heijnen\}@tudelft.nl). }
\thanks{P. Herder is with the Faculty of Applied Sciences, Delft University of Technology, Delft, the Netherlands (e-mail: p.m.herder@tudelft.nl). }
\thanks{This research received funding from the Netherlands Organisation for Scientific Research (NWO) [project number: 647.002.007].}

}

\maketitle

\begin{abstract}
The growing share of proactive actors in the electricity markets calls for more attention on prosumers and more support for their decision-making under decentralized electricity markets. In view of the changing paradigm, it is crucial to study the long-term planning under the decentralized and prosumer-centric markets to unravel the effects of such markets on the planning decisions.  
In the first part of the two-part paper, we propose a prosumer-centric framework for concurrent generation and transmission planning. Here, three planning models are presented where a peer-to-peer market with product differentiation, a pool market and a mixed bilateral/pool market and their associated trading costs are explicitly modeled, respectively. To fully reveal the individual costs and benefits, we start by formulating the optimization problems of various actors, i.e. prosumers, transmission system operator, energy market operator and carbon market operator. Moreover, to enable decentralized planning where the privacy of the prosumers is preserved, distributed optimization algorithms are presented based on the corresponding centralized optimization problems.  
\end{abstract}


\IEEEpeerreviewmaketitle

\section{Introduction}

To mitigate climate change and reduce carbon emissions, renewable energy sources (RES) play a vital role in modern power systems. Countries around the world set ambitious RES targets for the next decades and the planning of RES generation and transmission is prominent on the agenda. 

Integrated planning of generation and transmission assumes there is a centralized planning entity, who determines the optimal expansion of the system to fulfill the carbon target. This method serves as a benchmark, which, however, does not provide information on how it can be used in a market environment. In liberalized electricity markets, the generation investments are driven by the price signals from the markets.
It is thus essential to investigate the influences of the electricity markets on the planning decisions. 
In fact, the costs and benefits that are incurred from the electricity markets (later all referred to as trading-related costs, which can be negative in case of revenues), i.e. those associated with buying or selling energy in the markets, constitute a major part of the costs of the generation companies \cite{Botterud2005}, or in general the investors. They play a crucial role in their decision-making on whether or to what extent they will invest. 
However, in energy system models for integrated planning of generation and transmission, the markets are not yet modeled.

The other critical challenge in the prevailing centralized planning models is the lack of attention on the coordination problem when multiple areas need to plan together \cite{Wang2017}. Due to political, technical and privacy issues, it is difficult for the involving agents, e.g., the areas (or otherwise referred to as prosumers\footnote{In this study, the notation prosumer represents an agent with an aggregated supply and/or a demand, which is interchangeable with peer. It can refer to an agent on either the distribution level or the transmission level, which includes but is not limited to households.}) and the grid operators, to share data.  
With the large-scale integration of RES, the inter-dependencies that are enforced by interconnections between the areas are growing and thus proper coordination between the agents is called for to achieve overall economic efficiency. Therefore, a decentralized planning approach where the agents make decisions locally needs to be developed. This problem becomes more urgent, due to the emergence of more decentralized and prosumer-centric electricity markets, such as the peer-to-peer (P2P) electricity markets \cite{Parag2016}. 

Therefore, in this study, we propose a prosumer-centric framework for concurrent generation and transmission planning. Instead of the conventional integrated planning approach, we focus on the wills of the relevant agents, where they either do generation or transmission planning at the same time. This way of modeling is called a concurrent approach in this paper. To that end, various electricity markets are modeled explicitly in the planning models of the agents, including P2P markets, pool markets and mixed bilateral/pool markets. 
The proposed approach will help the prosumers and the grid operators to make informed planning decisions considering the electricity markets while accounting for the carbon target imposed by the government. It can also help policy-makers evaluate how various market designs would influence the decentralized planning decisions.

\subsection{Planning models integrating peer-to-peer electricity markets}

Traditional electricity markets are designed mainly to accommodate centralized generation sources. In recent years, with the increasing penetration of distributed energy sources, the market paradigms are changing and prosumer markets emerge as next-generation market designs. In these markets, prosumers actively participate in the market and are expected to gain substantial profits from the market. \cite{Parag2016} summarized three typical prosumer market models, i.e. P2P, prosumer-to-grid and organized prosumer groups. In particular, P2P markets have been studied extensively over the past few years.

A major focus of the existing studies of P2P electricity markets is on the operational perspective. 
The design of the market is recognized as one of the essential layers alongside the physical layer and information layer \cite{VanLeeuwen2020}. 
Here, bilateral trading is considered as one of the most promising P2P market mechanisms \cite{Wang2020c}. 
When prosumers trade bilaterally in P2P markets, they could express preferences for particular trades. There is not yet a naming convention for that in the field, terms such as heterogeneous preferences (\cite{Yang2015, Hahnel2020}), product differentiation \cite{Sorin2019}, and energy classes \cite{Morstyn2019a} have been used. Among those, product differentiation is a generic mathematical formulation \cite{Baroche2019b} that can be used for various purposes, e.g., \cite{Baroche2019a} used it to account for exogenous network charging in P2P markets.

Bilateral contracts exist not only in P2P markets, in fact, they constitute major parts of energy trading in the prevailing electricity markets. They are good at hedging against the uncertainties associated with the price and the quantity in day-ahead markets.
The bilateral contracts can be modeled directly by P2P market models. Besides, they are studied mainly using simulation methods (\cite{Bower2000, Knezevic2011, Lopes2012, Imran2020}) and complex network approaches \cite{Bompard2008}.

Yet, the investment decisions under the P2P markets are not discussed \cite{Wang2020c}. The same as the traditional pool electricity markets, the P2P markets also have a significant influence on the investment decisions of the prosumers by giving price signals. While the operational aspects of the P2P markets have been extensively studied, it is of importance to study how these markets would lead to what investment decisions. Therefore, a largely ignored aspect of the P2P markets, or in other words, a significant research gap is to integrate those markets in the planning models. Only by such models, prosumers can make informed planning decisions under a P2P market environment.

\subsection{Planning models integrating pool electricity markets}

Since the liberalization from the 1990s, the electricity sector has undergone drastic changes, where the pool-based electricity market has been the norm.
However, with the emergence of the decentralized electricity markets, the predominant paradigm has come to a turning point. Decentralized electricity markets are playing an increasingly crucial role, which might supersede the traditional market designs in the long term. While most studies investigate techno-economic uncertainties such as changes in demand, weather forecasts and cost parameters, the institutional uncertainty, i.e. the uncertainty in how the future electricity markets will be shaped, is generally overlooked. 
Planning models should highlight this uncertainty, unraveling the associated trading costs and benefits along the planning horizon. 
Consequently, it is of high importance to consider pool markets as well in the planning models, in addition to the P2P markets. The inclusion of both markets in one framework will show quantitatively and comparatively the incurred costs when different markets are considered. 
As an example, \cite{Doorman2008} analyzed the effects of four market designs on generation investment. 

In the literature, the decentralized planning of generation and transmission in the market environment is often addressed by multi-level optimization problems. From the game-theoretical perspective, these problems are also known as leader-follower (Stackelberg) games. In the upper level, the transmission system operator (TSO), as a leader, makes transmission investments, while anticipating the investments of generation companies. In the lower level, generation companies, as followers, maximize their profits, while their optimization problems are constrained by market-clearing outcomes. Several examples are \cite{Jenabi2013, Taheri2017, Grimm2020}. This type of studies takes the perspectives of the TSO and investigates the anticipatory sequential decision-making process. Similarly, some studies assume a reverse sequence in decision-making, i.e., the generation companies make decisions first taking into account the transmission charges, and then the TSO follows, e.g., see the work of \cite{Tohidi2017}. 

Instead of this sequential decision-making approach, the prosumers and the TSO may also make decisions at the same time, i.e. concurrently.
In P2P markets, the market-clearing is modeled as an equilibrium problem (\cite{Baroche2019b, Moret2020}), where all agents, i.e. prosumers, market operators and TSO, negotiate together while minimizing their costs. This approach is also used to study pool markets in \cite{Hobbs03complementarity-basedequilibrium, StevenA.Gabriel2013, Ruiz2014}. Moreover, it has also been applied to the planning problem under the pool market in \cite{He2012}.

\subsection{Distributed planning algorithms}

Unlike centralized optimization where information has to be gathered by the centralized entity, distributed optimization indicates that agents exchange information with each other and optimize their local problems in order to solve the centralized problem often without the presence of a centralized coordinator. \cite{Wang2017, Molzahn2017, Kargarian2018} give reviews on the distributed optimization approaches that are applied for power systems. Among those, Lagrangian relaxation and Benders decomposition together with their variants are the most commonly seen in the literature \cite{Sagastizabal2012}.

In planning models, Benders decomposition and its variants are often used (such as \cite{Ruiz2015,  Munoz2016}) to increase solution speed. 
In P2P market models, Lagrange relaxation, in particular, alternating direction method of multipliers (ADMM) has been used extensively. ADMM allows the agents to solve their local optimization problems while preserving private information. 
Distributed optimization is particularly needed in the P2P markets where preferences of the prosumers are considered, see e.g. \cite{Morstyn2019a, Sorin2019, Moret2019, Baroche2019b, LeCadre2020}. Given different application domains of the methods, the choice of a method depends on the purpose of using distributed optimization. In this study, since we aim to let the prosumers solve their local problems, ADMM will be further used and discussed.

\subsection{Contributions}
\begin{figure}[htbp]
\centerline{\includegraphics[width=0.5\textwidth]{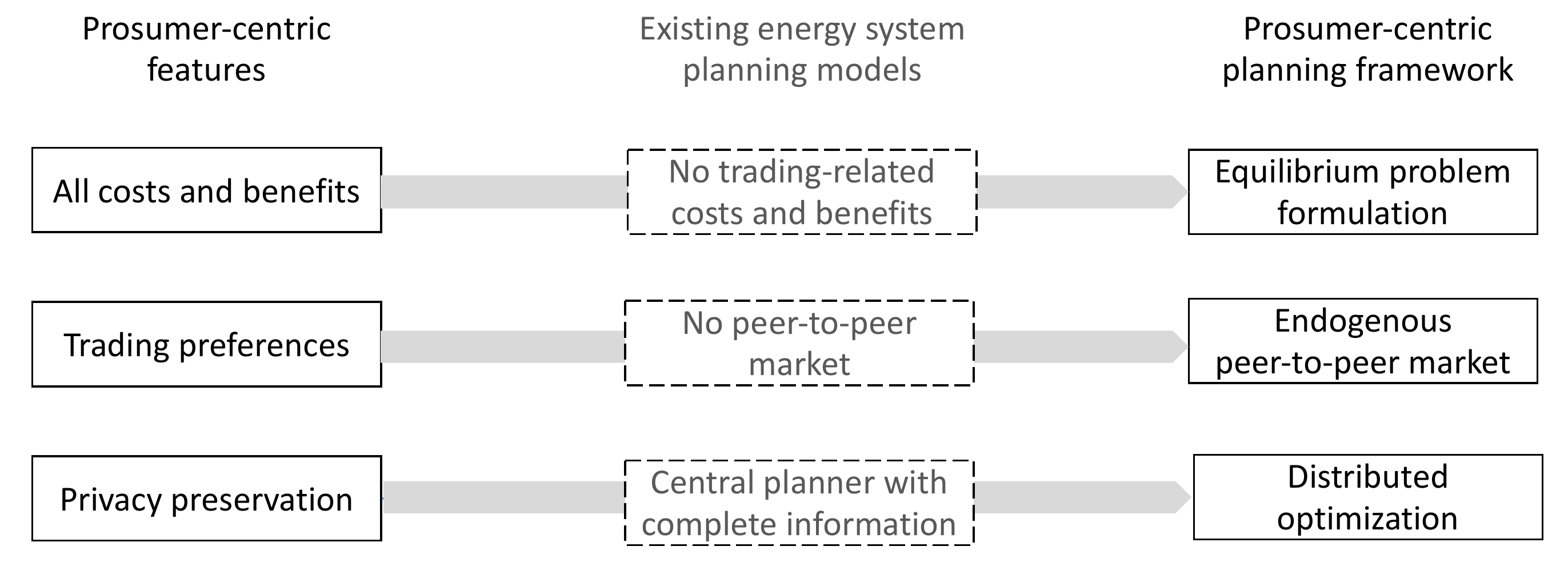}}
\caption{Prosumer-centric features and the proposed planning framework, compared to existing energy system planning models.}
\label{fig:features}
\end{figure}

Based on the background information and the literature review, we found that the existing energy system planning models have a focus on the system level. However, for modern power systems, especially with the increasing interests in P2P markets, the existing integrated planning models are not sufficient to capture the prosumer-level results and are thus not able to aid the investment decision-making of the prosumers. Therefore, we propose a prosumer-centric planning approach that incorporates the following features, which are also presented in Figure \ref{fig:features}.

\begin{enumerate}
    \item all costs and benefits, in particular, those associated with buying/selling energy on the electricity market are included in the objective functions of the prosumers. This requires the explicit formulation of the optimization problem for each agent. The problems are interconnected, and thus cannot be solved separately. To solve them together, an equilibrium problem has to be formulated.
    \item preferences of the prosumers are included. Compared to the objective where only costs are included, this results in an improved utility function that considers the willingness to pay of the prosumers. This is done by utilizing the production differentiation term from a P2P market. Consequently, the P2P market has to be considered endogenously in a planning model. 
    \item privacy of the prosumers is preserved. To fully account for the privacy of the prosumers, distributed optimization algorithms must be provided where each agent solves its problem locally with limited information exchange with each other. 
\end{enumerate}

In our proposed framework, besides the prosumers, the energy market operator and the TSO, the government is also included. The role of the government is to set carbon goals for the future, which are facilitated through policies such as a cap-and-trade system and carbon tax \cite{He2012}. In this study, a cap-and-trade system with a carbon market is considered since it is a market-based policy and a carbon target can be imposed directly by the government.

This planning framework includes three planning models which consider a P2P market, a pool market and a mixed bilateral/pool market, respectively. 
Part I of the paper focuses on model formulations, and numerical results will be provided in Part II.
The contributions of this paper are:
\begin{itemize}
    \item the formulation of a planning model integrating the P2P market design with product differentiation, where previously only operational problems are studied. This planning model is the first-of-its-kind that studies the P2P market beyond only the operational aspects.
    Such a model allows investigating the influence of the P2P market design on the planning decisions. 
    \item the formulation of planning models integrating the pool electricity market and the mixed bilateral/pool market, respectively.
    \item the formulations of the distributed optimization algorithms to solve the planning problems in a decentralized manner.  
\end{itemize}

The most common usage for energy system planning models is to assess the investment decisions and the associated costs. Naturally, the fundamental usage of our framework is to investigate the optimal investment decisions and the associated costs and benefits given the various market environments. Despite this, our framework is able to go beyond it and deal with other problems. Here, we will briefly introduce three potential applications and they will be further explained and illustrated by the archetypal case studies in Part II of the paper.

As a first example, in the modern and evolving power sector, the planning decisions have to be made ahead of time but there will be uncertainties in future market designs. Hence, it is important to make no-regret decisions such that the total cost is minimal. Thanks to the inclusion of various market designs, the framework is able to help the agents to best deal with this uncertainty.
Moreover, the framework can be used as a negotiation tool in joint generation planning of the prosumers. In the endogenous P2P market with product differentiation, the traded energy could be influenced by the willingness to pay of the prosumers, and thus affects the planning decisions. Here, the willingness to pay will play a crucial role in the investment negotiations. 
Furthermore, the framework allows to model bilateral contracts in the mixed bilateral/pool markets. 
It opens up research directions such as modeling the effects of bilateral contracts and investigating how product differentiation therein could influence the planning decisions.

The remainder of the study is structured as follows.
Section \ref{sec:pre} gives the preliminaries of the framework. Then, Section \ref{sec:p2p} - \ref{sec:mix} discuss the planning models integrating the P2P market, the pool market and the mixed bilateral/pool market, respectively. Lastly, conclusions are drawn in Section \ref{sec:conclusion}.

\section{Preliminaries} \label{sec:pre}

This section describes the preliminaries needed to understand the model and outlines
the objective functions for the prosumers and for the TSO in integrated energy system planning models. 

\begin{figure}[htbp]
\centerline{\includegraphics[width=0.5\textwidth]{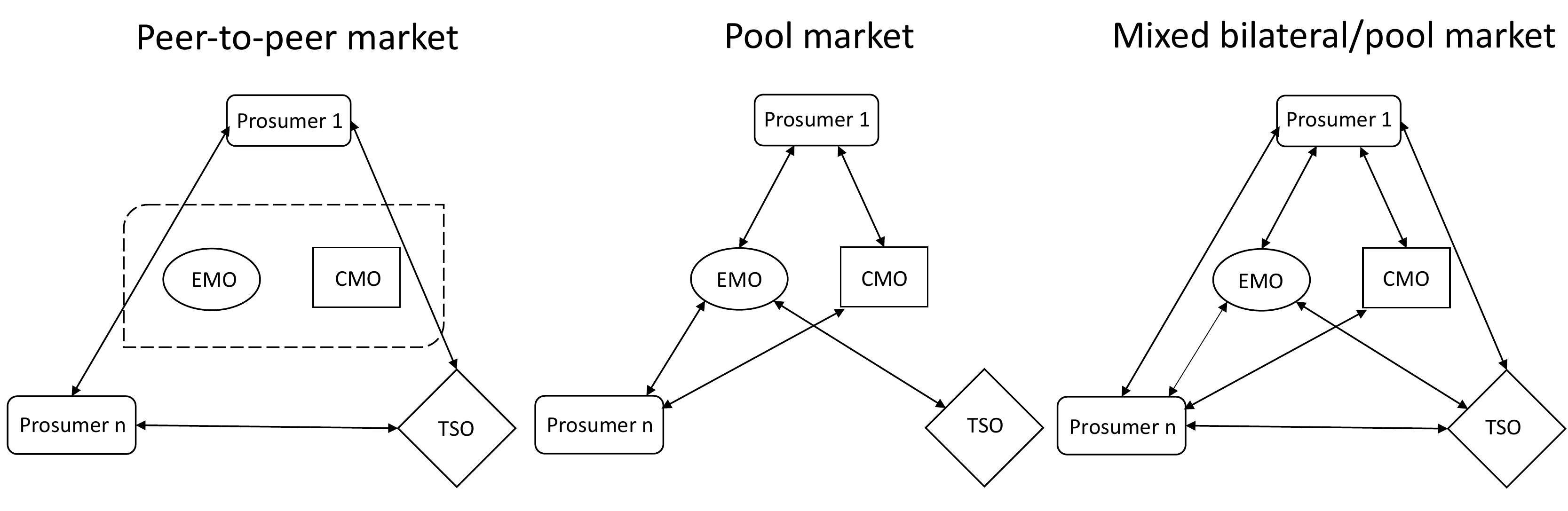}}
\caption{Information flows between the agents for peer-to-peer market, pool market and mixed bilateral/pool market. EMO: energy market operator, CMO: carbon market operator, TSO: transmission system operator. Dashed box indicates that EMO and CMO may not exist physically.}
\label{fig:structure}
\end{figure}

Let us start with briefly introducing the three market designs. Figure \ref{fig:structure} shows the structure of them in terms of the information flows.  
The P2P market electricity design is based on \cite{Baroche2019b, Moret2020}. In this P2P market, the prosumers are represented by the nodes in the communication graph and the edges connect the pair of the prosumers who might trade energy with each other. The prosumers can trade freely with each other that are on the same communication graph. The bilateral prices are based on negotiation and are derived using the proposed distributed algorithm which performs equally as if there are market operators physically present. Here, differently from the multi-level optimization that is commonly seen in the literature, the TSO is on the same level (and thus on the same communication graph) as the prosumers. In the conventional pool market, the market operators are physically present, and the information flows are only between them and the market participants. In the mixed bilateral/pool market, the market operators are in charge of the pool market while prosumers can also trade bilaterally under the negotiated prices. 

For each of the three markets, we start by formulating the optimization problems by assuming perfect competition for all the agents, i.e. the prosumers, the TSO, the energy market operator and the carbon market operator. This step is crucial because the individual costs can only be derived if the objective functions of the agents are formulated explicitly. Note that since these problems are all interconnected, i.e. the parameters in one problem may be the decision variables in the other, and vice versa, they should be solved together. This will be done by finding the Nash equilibrium, based on which the equivalent centralized optimization problems will be given. Lastly, distributed optimization algorithms are presented using ADMM.   

Now the terminology used in the paper will be introduced. Lower case symbols are used for variables, upper case symbols are used for parameters. Dual variables are expressed using Greek letters. $n$ is the index for prosumers $\mathcal{N}$. $i$ is the index for generation technologies $\mathcal{G}$, and storage technologies $\mathcal{S}$. $l$ is the index for transmission lines in the existing line set $\mathcal{L}$. $t$ represents a time step in $\mathcal{T}$. 

In the family of energy system planning models, the objective is usually to minimize total annualized cost consisting of Capital Expenditure (CapEx) cost of generation and storage technologies, Fixed Operation \& Maintenance (FOM) costs, Variable Operation \& Maintenance (VOM) costs and CapEx cost of networks \cite{Wang2020b}. 
These planning models use a direct split of the total costs that are incurred to the prosumers and the TSO. For the prosumer $n$, its cost $f_n$ is the summation of CapEx, FOM and VOM of the corresponding technologies. Here, $A_i$ is the annual factor for technology $i$, $C_i$ and $CS_i$ are the CapEx for generation $i$ and storage $i$, respectively, $k_{i,n}$ and $k^{\text{storage}}_{i,n}$ are the capacities for generation and storage technologies $i$ for prosumer $n$, respectively, $B_i$ is the VOM cost for technology $i$ and $p_{i,n,t}$ is the production of technology $i$ for prosumer $n$ at time step $t$. The TSO, however, has to bear the cost $g$, which is the total CapEx for the transmission network, where $\Delta_l$ is the length of line $l$.

\begin{subequations}
\begin{flalign}
	f_n & = \sum_{i\in (\mathcal{G}+\mathcal{S})} \frac{C_{i} k_{i,n}}{A_{i}} + \sum_{i\in \mathcal{S}} \frac{CS_{i} k^{\text{storage}}_{i,n}}{A_{i}} \\ 
	& + \sum_{t \in \mathcal{T}} \sum_{i\in \mathcal{G}}  B_i p_{i,n,t}  \notag \\
    g & = \sum_{l \in \mathcal{L}} \frac{\Delta_l C_l k_l}{A_l}
\end{flalign}
\end{subequations}

\section{Peer-to-peer market} \label{sec:p2p}
Under a P2P market, we start by formulating the interconnected optimization problems for all the agents. Next, the equivalent centralized optimization problem will be given. Then, the ADMM algorithm for distributed optimization will be described.   

\subsection{Prosumer}
{\allowdisplaybreaks

In this study, a prosumer is considered as an archetypal agent in the energy system with some generation capacities and/or an energy demand who participates in a P2P market on the transmission level. Examples can be groups of aggregators/large consumers/generation companies that are aggregated to a power injection point on the transmission network.
Note that, in principle, different prosumers can co-exist at a transmission power injection point, i.e. there might be a discrepancy between the set of prosumers and the set of transmission nodes. For simplicity but without loss of generality, we consider the two sets are identical. Moreover, a prosumer can be on a distribution level, e.g., a household and the gist of the framework still applies as long as TSO is changed to DSO, the market is changed to a local market and the power flow calculation is changed accordingly.

The prosumer $n \in \mathcal{N}$ can trade energy $p^{\text{p2p}}_{n,m,t} $ at time step $t$ bilaterally with its neighbors $m \in \omega_n$ that are on the communication graph. 
Trading prices $\lambda^{\text{p2p}}_{n,m,t}$ and grid prices $\lambda^{\text{grid}}_{n,m,t}$ are associated with the energy trades $p^{\text{p2p}}_{n,m,t} $. In addition, $I_{n,m}$ is the product differentiation parameter from $n$ to $m$. More details of the product differentiation term can be found in \cite{Sorin2019}.

$\Gamma^{\text{ps}}_n$ is the set of decision variables for prosumer $n$.
It includes the investment capacities $k_{i,n}$ of generation and storage conversion $i$, the investment capacities $k^{\text{storage}}_{i,n}$ of storage $i $, the energy production $p_{i,n,t}$ from technology $i$ at time step $t$, the  bilateral trades $p^{\text{p2p}}_{n,m,t} $ from $n$ to $m$ at time step $t$, state-of-charge $soc_{i,n,t}$ of storage $i$ at time step $t$, storage discharging $p^{\text{out}}_{i,n,t}$ of storage $i$ at time step $t$, storage charging $p^{\text{in}}_{i,n,t}$ of storage $i$ at time step $t$ and the number of carbon permits $e^{CO_2}_n$ to buy from the carbon market.

The objective function of prosumer $n$ is given in (\ref{equ:pro_obj}), which minimizes the total cost of the prosumer. The total cost consists of three parts. The first part $f_n$ is the cost that is commonly used in energy system planning models, i.e. CapEx, FOM and VOM. The second part is the trading-related costs in the energy market, including trading costs, grid costs and product differentiation costs. The third part is the trading-related costs in the carbon market. Since a cap-and-trade system is discussed here, the prosumer $n$ has to buy a certain amount of emissions permits $e^{CO_2}_n$ that are equivalent to their emissions. 
Therefore, the objective function (\ref{equ:pro_obj}) not only considers the planning costs but also integrates the trading-related costs (which can be negative provided that the prosumer gains more from selling energy) in the energy market and the carbon market. 

\begin{subequations}
\begin{flalign}
	\min_{\Gamma^{\text{ps}}_n}
	 f_n + \sum_{t \in \mathcal{T}} \sum_{m \in \omega_n} (\lambda^{\text{p2p}}_{n,m,t} + \lambda^{\text{grid}}_{n,m,t})  p^{\text{p2p}}_{n,m,t} \label{equ:pro_obj} \\
	 + \sum_{t \in \mathcal{T}} \sum_{m \in \omega_n} I_{n,m} |p^{\text{p2p}}_{n,m,t} | + \lambda^{CO_2} e^{CO_2}_n \notag 
\end{flalign}

\begin{flalign}
 	 \text{s.t.} & \sum_{i\in \mathcal{G}} p_{i,n,t} - D_{n, t} + \sum_{i\in \mathcal{S}}(p^{\text{out}}_{i,n,t} - p^{\text{in}}_{i,n,t})   \notag \\
 	 & = \sum_{m \in \omega_n} p^{\text{p2p}}_{n,m,t} ,  \forall t \in \mathcal{T}  \label{equ:pro_balance} \\
	 & 0 \le p_{i,n,t} \le E_{i,n,t} (k_{i,n} + K_{i,n})  , \forall i \in \mathcal{G}, \forall t \in \mathcal{T} \label{equ:pro_gen_limit} \\
	 & soc_{i,n,t} = soc_{i,n,t-1} + H^{in}_{i} p^{\text{in}}_{i,n,t} - \frac{1}{H^{out}_{i}} p^{\text{out}}_{i,n,t},  \notag\\
	 & \forall i \in \mathcal{S},  \forall t \in \mathcal{T} \label{equ:pro_storage_begin} \\
	 & 0 \le soc_{i,n,t} \le k^{\text{storage}}_{i,n} + K^{\text{storage}}_{i,n} , \forall i \in \mathcal{S}, \forall t \in \mathcal{T} \label{equ:cm4} \\
	 & 0 \le p^{\text{out}}_{i,n,t} \le k_{i,n} + K_{i,n}, \forall i \in \mathcal{S},  \forall t \in \mathcal{T} \\
	 & 0 \le p^{\text{in}}_{i,n,t} \le k_{i,n} + K_{i,n}, \forall i \in \mathcal{S},  \forall t \in \mathcal{T} \label{equ:pro_storage_end} \\
	 & e^{CO_2}_n = W_i \sum_{t \in \mathcal{T}} \sum_{i \in \text{$\mathcal{R}$}} p_{i,n,t}  \label{equ:pro_CO2}
\end{flalign}
\end{subequations}

}

(\ref{equ:pro_balance}) is the nodal energy balance constraint. On the left-hand side of the equation is the net energy that can be traded, and on the right-hand is the sum of the bilateral trades.
(\ref{equ:pro_gen_limit}) indicates that the energy production is constrained by the efficiency $E_{i,n,t}$ (capacity factor in case of variable renewable energy) and the capacity of the generation technologies. Here, $K_{i,n}$ is the existing capacity, $k_{i,n}$ is the capacity to be expanded, essentially making the model a generation expansion model. (\ref{equ:pro_storage_begin}) - (\ref{equ:pro_storage_end}) are the storage constraints, as described in \cite{Wang2020b}. (\ref{equ:pro_CO2}) shows that the amount of emissions is equal to the required carbon permits.

\subsection{TSO}

The role of the TSO is two-fold. On the one hand, it ensures the feasibility of the energy flows and accordingly, invests in the transmission network capacity cost-optimally. 
On the other hand, it harvests congestion rents by trading $z^{\text{p2p}}_{n,m,t}$ on the electricity market as a spatial arbitrager. As a result, the objective function (\ref{equ:tso_obj}) is to minimize the total cost which incorporates the cost pertaining to these two roles respectively.
The decision variables of TSO are represented by the set $\Gamma^{\text{TSO}}$, which includes the investment capacity $k_l$ in line $l$, the trades $z^{\text{p2p}}_{n,m,t}$ from $n$ to $m$ at times step $t$ and the energy flow $f_{l,t}$ in line $l$ at times step $t$.

\begin{subequations}
\begin{flalign}
	\min_{\Gamma^{\text{TSO}}} g - \sum_{t \in \mathcal{T}} \sum_{n \in \mathcal{N}} \sum_{m \in \omega_ n} \lambda^{\text{grid}}_{n,m,t} z^{\text{p2p}}_{n,m,t} \label{equ:tso_obj}
\end{flalign}
\begin{flalign}
     \text{s.t.} & f_{l,t} = \sum_{n \in \mathcal{N}} PTDF_{l,n} \sum_{m \in \omega_n} z^{\text{p2p}}_{n,m,t},  \forall l \in \mathcal{L},  \forall t \in \mathcal{T} \label{equ:tso_ptdf} \\
	 & - (k_l + K_l) \le f_{l,t} \le k_l + K_l,  \forall l \in \mathcal{L},  \forall t \in \mathcal{T} \label{equ:tso_limit}
\end{flalign}
\end{subequations}

The energy flow is modeled using DC power flow equations. In (\ref{equ:tso_ptdf}), the flow $f_{l,t}$ is calculated based on the Power Transfer Distribution Factors (PTDF) matrix and the total net injection at $n$. Except for this formulation, other formulations of DC power flow can also be used. Moreover, AC power flows might be considered as well, e.g., \cite{Moret2020} studied the network losses in a joint transmission and distribution P2P markets, which can be used to complement our framework to consider network operators both on the transmission and distribution level. 
(\ref{equ:tso_limit}) indicates the thermal limits of the energy flows, where $k_l$ refers to the transmission capacity expansion. 

\subsection{P2P market operator}

The P2P energy market operator clears the market at each time step $t$ by minimizing the energy imbalances, and thus determines the corresponding prices. The set of decision variables $\Gamma^{\text{market}}$ includes the trading price $\lambda^{\text{p2p}}_{n,m,t}$ from $n$ to $m$ at time step $t$ and the grid price $\lambda^{\text{grid}}_{n,m,t}$ from $n$ to $m$ at time step $t$. It makes sure that the bilateral trades should be equal in quantity, and the trading energy from the prosumer $p^{\text{p2p}}_{n,m,t} $ is equal to the arbitraging energy from the TSO $z^{\text{p2p}}_{n,m,t}$ at each time step $t$.

\begin{flalign}
	\min_{\Gamma^{\text{market}}} \sum_{t \in \mathcal{T}} \sum_{n \in \mathcal{N}} \sum_{m \in \omega_n} \lambda^{\text{p2p}}_{n,m,t} (p^{\text{p2p}}_{n,m,t} +  p^{\text{p2p}}_{m,n,t}) \notag \\
	+ \sum_{t \in \mathcal{T}} \sum_{n \in \mathcal{N}} \sum_{m \in \omega_n} \lambda^{\text{grid}}_{n,m,t} (p^{\text{p2p}}_{n,m,t} -z^{\text{p2p}}_{n,m,t})
\end{flalign}

\subsection{Carbon market operator}
We consider a cap-and-trade system as the carbon policy. The government can determine the maximum amount of emissions that are allowed to be emitted, which will be regarded as a cap on the carbon market. All prosumers need to buy carbon permits that are equivalent to their emissions from the carbon market. The price $\lambda^{CO_2}$ will be determined by the carbon market. 
It is formulated as the following. 

\begin{flalign}
	\min_{\lambda^{CO_2}} \lambda^{CO_2} (\sum_{n \in \mathcal{N}} e^{CO_2}_n - \text{CAP}^{CO_2})
\end{flalign}

\subsection{Equivalent planning optimization problem}

After laying down the optimization problems of the prosumers, the TSO, the energy market operator and the carbon market operator, it could be found that using this way of formulation, the decision variables of one problem only exist in the objective function of others and not in the constraints which avoids resulting in a generalized Nash equilibrium. We are then able to find the Nash equilibrium by deriving the the Karush–Kuhn–Tucker (KKT) conditions of these problems. The equivalent optimization problem can be written as follows: 

{\allowdisplaybreaks
\begin{subequations}
\begin{flalign}
	\min_{\Gamma}
	 \sum_{n \in \mathcal{N}} f_n +  g + \sum_{t \in \mathcal{T}} \sum_{n \in \mathcal{N}} \sum_{m \in \omega_n} I_{n,m} |p^{\text{p2p}}_{n,m,t} | \label{equ:p2p_cen_obj}
\end{flalign}
\begin{flalign}
    \text{s.t.}  & (\ref{equ:pro_balance}) - (\ref{equ:pro_CO2}), \forall n \in \mathcal{N}\\
     & (\ref{equ:tso_ptdf}) - (\ref{equ:tso_limit}) \\
	 & p^{\text{p2p}}_{n,m,t} = - p^{\text{p2p}}_{m,n,t}, \forall n \in \mathcal{N}, \forall m \in \omega_n, \forall t\in T  \text{: } \lambda^{\text{p2p}}_{n,m,t} \label{equ:mo_reciprocity} \\
	 & p^{\text{p2p}}_{n,m,t} = z^{\text{p2p}}_{n,m,t},  \forall n \in \mathcal{N}, \forall m \in \omega_n,\forall t\in T   \text{: } \lambda^{\text{grid}}_{n,m,t}  \label{equ:mo_grid} \\
	 & \sum_{n \in \mathcal{N}} e^{CO_2}_n = \text{CAP}^{CO_2} \text{: } \lambda^{CO_2} \label{equ:mo_co2}
\end{flalign}
\end{subequations}
}

The decision variables belong to the set $\Gamma$, which includes all the decision variables of the prosumers and the TSO. 
The objective function (\ref{equ:p2p_cen_obj}) is the summation of the objective functions of all the agents. 
The constraints are also a gathering of all constraints of the agents' problems. Note that (\ref{equ:mo_reciprocity}) and (\ref{equ:mo_grid}) are KKT conditions of the optimization problem of the energy market operator. (\ref{equ:mo_reciprocity}) is the reciprocity constraint, showing that the bilateral trades should be equal in quantity, where the dual variable $\lambda^{\text{p2p}}_{n,m,t}$ is the trading price. (\ref{equ:mo_grid}) is the energy balance constraint at $n$, equalizing the bilateral trades from prosumer $p^{\text{p2p}}_{n,m,t} $ to the arbitraging energy from the TSO $z^{\text{p2p}}_{n,m,t}$, where the dual variable is the grid price for this trade. (\ref{equ:mo_co2}) gives the cap for all the carbon emissions, with the dual variable $\lambda^{CO_2}$ being the carbon price.

\subsection{Distributed solving formulation using ADMM}


ADMM is a distributed optimization algorithm that gains popularity in recent years, for its simple formulation and guaranteed convergence. It is widely used as a solving technique in P2P markets. Although different algorithms can be used to solve the planning optimization problem (\ref{equ:p2p_cen_obj}), we chose here to employ ADMM because it fits better in the proposed prosumer-centric framework as it allows decentralized decision-making.

The Lagrangian of (\ref{equ:p2p_cen_obj}) is formulated in (\ref{equ:p2p_admm_L}). Compared to (\ref{equ:p2p_cen_obj}), the coupling constraints (\ref{equ:mo_reciprocity})-(\ref{equ:mo_co2}), i.e. the constraints that contain decisions variables of more than one agent, are relaxed. Three quadratic terms with the penalty parameters $Q^{\text{trade}}_{n,m,t}$, $Q^{\text{grid}}_{n,m,t}$ and $Q^{CO_2}$ are added to ensure convergence. 

{\allowdisplaybreaks
\begin{subequations}
\begin{flalign}
    \mathcal{L}(\Gamma^{\text{ps}}_n, \Gamma^{\text{TSO}}) = \sum_{n \in \mathcal{N}} f_n + g \notag 
    + \sum_{t \in \mathcal{T}} \sum_{n \in \mathcal{N}} \sum_{m \in \omega_n} I_{n,m} |p^{\text{p2p}}_{n,m,t} | \notag \\
     +	\sum_{t \in \mathcal{T}} \sum_{n \in \mathcal{N}} \sum_{m \in \omega_n} \lambda^{\text{p2p}}_{n,m,t}  (p^{\text{p2p}}_{n,m,t} + p^{\text{p2p}}_{m,n,t}) \notag \\
     +	\sum_{t \in \mathcal{T}} \sum_{n \in \mathcal{N}} \sum_{m \in \omega_n} \frac{Q^{\text{p2p}}_{n,m,t}}{2} (p^{\text{p2p}}_{n,m,t} + p^{\text{p2p}}_{m,n,t})^2 \notag \\
     + \sum_{t \in \mathcal{T}} \sum_{n \in \mathcal{N}} \sum_{m \in \omega_n} \lambda^{\text{grid}}_{n,m,t} (p^{\text{p2p}}_{n,m,t} - z^{\text{p2p}}_{n,m,t}) \notag \\
     + \sum_{t \in \mathcal{T}} \sum_{n \in \mathcal{N}} \sum_{m \in \omega_n} \frac{Q^{\text{grid}}_{n,m,t}}{2} (p^{\text{p2p}}_{n,m,t} - z^{\text{p2p}}_{n,m,t})^2 \notag \\
     +	\lambda^{CO_2} (\sum_{n \in \mathcal{N}} e^{CO_2}_n - \text{CAP}^{CO_2}) \notag \\
     +	\frac{Q^{CO_2}}{2}  (\sum_{n \in \mathcal{N}} e^{CO_2}_n - \text{CAP}^{CO_2})^2 \label{equ:p2p_admm_L}
\end{flalign}

Then the ADMM algorithm is formulated in the following way that the decision-making of the four types of agents, i.e. the prosumers, the TSO and the energy market operator and the carbon market operator are all done locally. More specifically, in (\ref{equ:p2p_admm_pro}) and (\ref{equ:p2p_admm_tso}), the prosumers and the TSO solve their local optimization problems, respectively. Then the energy market operator updates the trade price and the grid price in (\ref{equ:p2p_admm_trade}) and (\ref{equ:p2p_admm_grid}), and the carbon market operator updates the carbon price in (\ref{equ:p2p_admm_co2}). The problems will be solved in an iterative manner, where $k$ is the number of iterations.

\begin{flalign}
    \Gamma^{\text{ps}(k+1)}_n = \text{argmin}_{\Gamma^{\text{ps}}_n} \mathcal{L}(\Gamma^{\text{ps}}_n, \Gamma^{\text{TSO}(k)}), \forall n \in \mathcal{N}  \label{equ:p2p_admm_pro} \\
    \text{s.t. } (\ref{equ:pro_balance}) - (\ref{equ:pro_CO2}) \notag \\    
    \Gamma^{\text{TSO}(k+1)} = \text{argmin}_{\Gamma^{\text{TSO}}} \mathcal{L}(\Gamma^{\text{ps}(k)}_n, \Gamma^{\text{TSO}}) \label{equ:p2p_admm_tso}\\
    \text{s.t. } (\ref{equ:tso_ptdf}) - (\ref{equ:tso_limit}) \notag \\
    \lambda^{\text{p2p}(k+1)}_{n,m,t} = \lambda^{\text{p2p}(k)}_{n,m,t} +  Q^{\text{trade}}_{n,m,t} (p^{\text{p2p}(k+1)}_{n,m,t} + p^{\text{p2p}(k+1)}_{m,n,t}),  \notag \label{equ:p2p_admm_trade}\\
    \forall t \in \mathcal{T}, \forall n \in \mathcal{N}, \forall m \in \omega_n \\
    \lambda^{\text{grid}(k+1)}_{n,m,t} = \lambda^{\text{grid}(k)}_{n,m,t} +  Q^{\text{grid}}_{n,m,t} (p^{\text{grid}(k+1)}_{n,m,t} - z^{\text{p2p}(k+1)}_{n,m,t}), \notag \label{equ:p2p_admm_grid} \\
    \forall t \in \mathcal{T}, \forall n \in \mathcal{N}, \forall m \in \omega_n \\
    \lambda^{{CO_2}(k+1)} = \lambda^{{CO_2}(k)} + Q^{CO_2}  (\sum_{n \in \mathcal{N}} e^{{CO_2}(k+1)}_n - \text{CAP}^{CO_2}) \label{equ:p2p_admm_co2}
\end{flalign}

\end{subequations}
}

\section{Pool market} \label{sec:pool}
Since the key of our contribution is to integrate different electricity market designs into the planning problems, in this section, we turn to another market design, i.e. the pool market. As we did for the P2P market, we will first formulate the optimization problems of the four types of agents, then give the equivalent centralized optimization problem and finally describe the ADMM algorithm for distributed optimization.

\subsection{Prosumer}

{\allowdisplaybreaks

In the pool market, the main difference with the P2P market is that there are no bilateral trades, instead, the energy trading will be done through a pool. 
The decision variables of prosumer $n$ in the pool market is similar to those in the P2P market, only without $p^{\text{p2p}}_{n,m,t}$.
The objective function of the prosumer $n$ is to minimize its total cost consisting of $f_n$ and the trading-related costs. The net power injection of prosumer $n$ is traded in the pool, under a nodal electricity price $\lambda^{\text{pool}}_{n,t}$ at time step $t$.

\begin{subequations}
\begin{flalign}
	\min_{\Gamma^{\text{ps}}_n}
	 f_n & + \sum_{t \in \mathcal{T}} \lambda^{\text{pool}}_{n,t} (\sum_{i\in \mathcal{G}} p_{i,n,t} - D_{n, t}+ \sum_{i\in \mathcal{S}}p^{\text{out}}_{i,n,t} - \sum_{i\in \mathcal{S}}p^{\text{in}}_{i,n,t}) \notag \\
	 & + \lambda^{CO_2} e^{CO_2}_n  \label{equ:nodal_pro_obj}
\end{flalign}

\begin{flalign}
 \text{s.t.}  (\ref{equ:pro_gen_limit}) - (\ref{equ:pro_CO2})
\end{flalign}
\end{subequations}

This optimization problem has the same constraints as the prosumer problem in the P2P market regarding generation limits, storage and carbon permits. However, the energy balance constraint is removed, since the net power injection is directly used in the objective function.

\subsection{TSO}

\begin{subequations}

The TSO has the same role as in the P2P market. It acts as a spatial arbitrager, ensures the feasibility of the energy flows and makes network investments accordingly. However, instead of $z^{\text{p2p}}_{n, m, t}$,
the TSO now only trades $z^{\text{pool}}_{n, t}$ at $n$ at time step $t$ at the price $\lambda^{\text{pool}}_{n,t}$. Its set of decision variables $\Gamma^{\text{TSO}}$ now includes the investment capacity $k_l$, the traded energy $z^{\text{pool}}_{n, t}$ and the energy flow $f_{l,t}$.

\begin{flalign}
	\min_{\Gamma^{\text{TSO}}} g - \sum_{t \in \mathcal{T}} \sum_{n \in \mathcal{N}}  \lambda^{\text{pool}}_{n,t}  z^{\text{pool}}_{n, t} \label{equ:nodal_tso_obj}
\end{flalign}

\begin{flalign}
     \text{s.t.} & \sum_{n \in \mathcal{N}} z^{\text{pool}}_{n, t} = 0, \forall t \in \mathcal{T} \label{equ:nodal_tso_balance} \\
     & f_{l,t} = \sum_{n \in \mathcal{N}} \text{PTDF}_{l,n} z^{\text{pool}}_{n, t} , \forall l \in \mathcal{L}, \forall t \in \mathcal{T} \label{equ:nodal_tso_flow}\\
	 & - (k_l + K_l) \le f_{l,t} \le k_l + K_l, \forall l \in \mathcal{L}, \forall t \in \mathcal{T} \label{equ:nodal_tso_limit} 
\end{flalign}
\end{subequations}

Since TSO is only a spatial arbitrager, (\ref{equ:nodal_tso_balance}) indicates that its total energy trades at $t$ should be equal to zero. (\ref{equ:nodal_tso_flow}) and (\ref{equ:nodal_tso_limit}) shows the DC power flow calculations using PTDF matrix and the thermal limits of the lines, respectively.

\subsection{Pool market operator}
The role of the energy market operator is to warrant the energy balance at each node $n$ and to derive the nodal price. Its set of decisions variables $\Gamma^{\text{pool}}$ now only contains $\lambda^{\text{pool}}_{n,t}$ which is the price for prosumer $n$ at time step $t$. 
It solves the following optimization problem.

\begin{flalign}
	\min_{\Gamma^{\text{pool}}}  \sum_{t \in \mathcal{T}} \sum_{n \in \mathcal{N}} \lambda^{\text{pool}}_{n,t} 
	& (\sum_{i\in \mathcal{G}} p_{i,n,t} - D_{n, t} + \sum_{i\in \mathcal{S}}p^{\text{out}}_{i,n,t} - \sum_{i\in \mathcal{S}}p^{\text{in}}_{i,n,t} \notag \\
	&- z^{\text{pool}}_{n, t}) \label{equ:nodal_market_obj}
\end{flalign}

\subsection{Equivalent planning optimization problems}


The carbon market operator's problem is identical to the problem under the P2P market design.
And similarly to deriving the equivalent planning optimization for the P2P market, the centralized optimization problem for the pool market can be derived as well. The solution to this problem provides the generation and transmission investment equilibrium under the pool electricity market. The objective function is a summation of (\ref{equ:nodal_pro_obj}), (\ref{equ:nodal_tso_obj}) and (\ref{equ:nodal_market_obj}). 
Note that this objective function is the same as the commonly-used objective functions in energy system planning models. In the literature, as we explained in Section \ref{sec:pre}, a direct split of the costs was often used. However, this objective function alone or this way of splitting costs do not reveal the real costs, i.e. those including the trading-related costs, for each agent. 
Only after we have shown the objective functions of each agent, it is now clear that their objectives sum up to this objective (with offsets in trading costs). Due to our prosumer-centric focus, we derived this relationship, such that the solutions will provide more insights from the agents' perspectives rather than only from the system perspective.    

\begin{subequations}
\begin{flalign}
	\min_{\Gamma}
	 \sum_{n \in \mathcal{N}} f_n + g \label{equ:pool_cen_obj}
\end{flalign}
\begin{flalign}
\text{s.t.}  & (\ref{equ:pro_balance}) - (\ref{equ:pro_CO2}), \forall n \in \mathcal{N} \label{equ:nodal_central_pro} \\
     & (\ref{equ:nodal_tso_balance}) - (\ref{equ:nodal_tso_limit}) \label{equ:nodal_central_tso}\\
	 & \sum_{i\in \mathcal{G}} p_{i,n,t} - D_{n, t} + \sum_{i\in \mathcal{S}}(p^{\text{out}}_{i,n,t} - p^{\text{in}}_{i,n,t})  = z^{\text{pool}}_{n, t}, \notag \\
	 & \forall n \in \mathcal{N}, \forall t \in \mathcal{T} \text{: }\lambda^{\text{pool}}_{n,t} \label{equ:nodal_central_balance} \\
	 & \sum_{n \in \mathcal{N}} e^{CO_2}_n = \text{CAP}^{CO_2} \text{: } \lambda^{CO_2} \label{equ:pool_cen_co2}
\end{flalign}
\end{subequations}

The set of decision variables of this problem $\Gamma$ includes the decision variables for all the prosumers and the TSO.
(\ref{equ:nodal_central_pro}) and (\ref{equ:nodal_central_tso}) are the constraints from the prosumers and the TSO, respectively. (\ref{equ:nodal_central_balance}) warrants that the energy is balanced for each node $n$ at time step $t$, i.e. the energy traded by the prosumer equals the arbitraged energy by the TSO. The dual variable of this constraint is the nodal energy price. Note here, different from the P2P market, there is no specific grid price to pay as it is part of the nodal price.

\subsection{Distributed solving formulation using ADMM}
The Lagrangian of the centralized planning problem (\ref{equ:pool_cen_obj}) is given in (\ref{equ:pool_admm_L}). The two coupling constraints (\ref{equ:nodal_central_balance}) and (\ref{equ:pool_cen_co2}) are relaxed, in which the penalty parameters $Q^{\text{pool}}_{n,t}, Q^{CO_2}$ are associated with the quadratic terms.


{\allowdisplaybreaks
\begin{subequations}
\begin{flalign}
    \mathcal{L}(\Gamma^{\text{ps}}_n, \Gamma^{\text{TSO}}) & = \sum_{n \in \mathcal{N}} f_n + g  \notag \\
    & +	\sum_{t \in \mathcal{T}} \sum_{n \in \mathcal{N}}  \lambda^{\text{pool}}_{n,t} (\sum_{i\in \mathcal{G}} p_{i,n,t} - D_{n, t} \notag \\
    & + \sum_{i\in \mathcal{S}}p^{\text{out}}_{i,n,t} - \sum_{i\in \mathcal{S}}p^{\text{in}}_{i,n,t} - z^{\text{pool}}_{n, t})  \notag \\
    & +	\sum_{t \in \mathcal{T}} \sum_{n \in \mathcal{N}}  \frac{Q^{\text{pool}}_{n,t}}{2}  (\sum_{i\in \mathcal{G}} p_{i,n,t} - D_{n, t} \notag \\
    & + \sum_{i\in \mathcal{S}}p^{\text{out}}_{i,n,t} - \sum_{i\in \mathcal{S}}p^{\text{in}}_{i,n,t}  - z^{\text{pool}}_{n, t})^2 \notag \\
    & +	\lambda^{CO_2} (\sum_{n \in \mathcal{N}} e^{CO_2}_n - \text{CAP}^{CO_2}) \notag \\ 
    & +	\frac{Q^{CO_2}}{2}  (\sum_{n \in \mathcal{N}} e^{CO_2}_n - \text{CAP}^{CO_2})^2 \label{equ:pool_admm_L}
\end{flalign}

The following ADMM algorithm features the local decision-making of the prosumers and the TSO, where their solutions will be sent to the market operators to update the price which will be given back to the prosumers and the TSO for iteration $k$ until convergence for all the problems.

\begin{flalign}
    \Gamma^{\text{ps}(k+1)}_n = \text{argmin}_{\Gamma^{\text{ps}}_n} \mathcal{L}(\Gamma^{\text{ps}}_n, \Gamma^{\text{TSO}(k)}), \forall n \in \mathcal{N} \\
    \text{s.t. } (\ref{equ:pro_gen_limit}) - (\ref{equ:pro_CO2}) \notag \\
    \Gamma^{\text{TSO}(k+1)} = \text{argmin}_{\Gamma^{\text{TSO}}} \mathcal{L}(\Gamma^{\text{ps}(k)}_n, \Gamma^{\text{TSO}})\\
    \text{s.t. } (\ref{equ:nodal_tso_balance}) - (\ref{equ:nodal_tso_limit}) \notag \\
    \lambda^{\text{pool}(k+1)}_{n,t} = \lambda^{\text{pool}(k)}_{n,t} + Q^{\text{pool}}_{n,t}  (\sum_{i\in \mathcal{G}} p_{i,n,t}^{(k+1)} - D_{n, t} + \sum_{i\in \mathcal{S}} p^{\text{out}(k+1)}_{i,n,t} \notag \\
    - \sum_{i\in \mathcal{S}}p^{\text{in}(k+1)}_{i,n,t} - z^{\text{pool}(k+1)}_{n, t}),  \forall n \in \mathcal{N},  \forall t \in \mathcal{T}  \\
    \lambda^{{CO_2}(k+1)} = \lambda^{{CO_2}(k)} + Q^{CO_2}  (\sum_{n \in \mathcal{N}} e^{{CO_2}(k+1)}_n - \text{CAP}^{CO_2}) 
\end{flalign}
\end{subequations}

}

\section{Mixed bilateral/pool market} \label{sec:mix}

In this section, we will integrate a mixed pool/bilateral electricity market into the planning model. At first, the optimization problems of all the agents are formulated. Then, an equivalent centralized optimization model is given. Finally, an ADMM algorithm is presented.

\subsection{Prosumer}
The prosumer $n$ aims to minimize the total costs that consist of $f_n$ and trading-related costs. Here, the trading-related costs in the energy markets include those in both the P2P market and the pool market. Consequently, the decision variables now include all those in the prosumer problems both markets. In particular, the prosumer $n$ can trade bilaterally which is represented by $p^{\text{bi}}_{n,m,t}$, but also in the pool by $p^{\text{pool}}_{n,t}$, representing the traded energy in the pool market for prosumer $n$ at time step $t$. 

\begin{subequations}
\begin{flalign}
	\min_{\Gamma^{\text{ps}}_n}
	 & f_n + \sum_{t \in \mathcal{T}} \sum_{m \in \omega_n}(\lambda^{\text{bi}}_{n,m,t} +  \lambda^{\text{grid}}_{n,m,t})  p^{\text{bi}}_{n,m,t} \notag \\
	 & + \sum_{t \in \mathcal{T}} \sum_{m \in \omega_n} I_{n,m} |p^{\text{bi}}_{n,m,t}|  
	  + \sum_{t \in \mathcal{T}} \lambda^{\text{pool}}_{n,t} p^{\text{pool}}_{n,t} + \lambda^{CO_2} e^{CO_2}_n
\end{flalign}  
\begin{flalign}
 \text{s.t.} & \Phi_n (\sum_{i\in \mathcal{G}} p_{i,n,t} - D_{n, t} + \sum_{i\in \mathcal{S}}p^{\text{out}}_{i,n,t} - \sum_{i\in \mathcal{S}}p^{\text{in}}_{i,n,t}) \notag \\
 & = \sum_{m \in \omega_n} p^{\text{bi}}_{n,m,t},  \forall t \in \mathcal{T}  \label{equ:mix_pro_balance1} \\
 	 & (1 - \Phi_n) (\sum_{i\in \mathcal{G}} p_{i,n,t} - D_{n, t} + \sum_{i\in \mathcal{S}}p^{\text{out}}_{i,n,t} - \sum_{i\in \mathcal{S}}p^{\text{in}}_{i,n,t}) \notag \\
 	 & = p^{\text{pool}}_{n,t},  \forall t \in \mathcal{T} \label{equ:mix_pro_balance2} \\  
 & (\ref{equ:pro_gen_limit}) - (\ref{equ:pro_CO2}) \label{equ:mix_pro_others}
\end{flalign}
\end{subequations}

(\ref{equ:mix_pro_balance1}) and (\ref{equ:mix_pro_balance2}) are both the energy balance constraints.  The net power injection $\sum_{i\in \mathcal{G}} p_{i,n,t} - D_{n, t} + \sum_{i\in \mathcal{S}}(p^{\text{out}}_{i,n,t} - p^{\text{in}}_{i,n,t})$ is divided into two parts: one part for the bilateral trading and  the other part for the pool-based trading. $\Phi_n$ is a parameter between 0 - 1 that is determined by the prosumer $n$ itself, indicating the percentage of its net energy that $n$ would like to trade bilaterally, the rest will be traded in the pool. In other words, the prosumer has to decide ex-ante how much to trade in the P2P market and in the pool market.
(\ref{equ:mix_pro_others}) gives the technical constraints on generation, storage and carbon permits. 

\subsection{TSO}
In addition to the investment cost $g$, the TSO receives the congestion rent from both electricity markets. The objective function is given in (\ref{equ:mix_tso_obj}). (\ref{equ:mix_tso_ptdf}) and (\ref{equ:mix_tso_limit}) ensure the feasibility of the flows. 

\begin{subequations}
\begin{flalign}
	\min_{\Gamma^{\text{TSO}}} g - \sum_{t \in \mathcal{T}} \sum_{n \in \mathcal{N}}( \sum_{m \in \omega_ n}  \lambda^{\text{grid}}_{n,m,t} z^{\text{bi}}_{n,m,t} -  \lambda^{\text{pool}}_{n,t}  z^{\text{pool}}_{n, t} ) \label{equ:mix_tso_obj}
\end{flalign}
\begin{flalign}
     \text{s.t.} & f_{l,t} = \sum_{n \in \mathcal{N}} PTDF_{l,n} (\sum_{m \in \omega_n} z^{\text{bi}}_{n,m,t} + z^{\text{pool}}_{n, t}),  \forall l \in \mathcal{L},  \forall t \in \mathcal{T} \label{equ:mix_tso_ptdf} \\
	 & - (k_l + K_l) \le f_{l,t} \le k_l + K_l,  \forall l \in \mathcal{L},  \forall t \in \mathcal{T} \label{equ:mix_tso_limit}
\end{flalign}
\end{subequations}

\subsection{Bilateral market operator}
The bilateral market operator ensures the energy balance for the bilateral trades and derives the associated energy and grid prices, by solving the following optimization problem. Its set of decision variables $\Gamma^{\text{bi}}$ consists of the trading price $\lambda^{\text{bi}}_{n,m,t}$ and the grid price $\lambda^{\text{grid}}_{n,m,t}$ with regarding to the bilateral trade from prosumer $n$ to prosumer $m$ at time step $t$.

\begin{flalign}
	\min_{\Gamma^{\text{bi}}} & \sum_{t \in \mathcal{T}} \sum_{n \in \mathcal{N}} \sum_{m \in \omega_n} \lambda^{\text{bi}}_{n,m,t} (p^{\text{bi}}_{n,m,t} +  p^{\text{bi}}_{m,n,t}) \notag \\
	&+ \sum_{t \in \mathcal{T}} \sum_{n \in \mathcal{N}} \sum_{m \in \omega_n} \lambda^{\text{grid}}_{n,m,t} (p^{\text{bi}}_{n,m,t} -z^{\text{bi}}_{n,m,t})   \label{equ:mix_market_obj1}
\end{flalign}

\subsection{Pool market operator}
The pool market operator has the same role as that of the bilateral market operator, but only for the energy that is traded in the pool. The following optimization problem will be solved. Its set of decision variables $\Gamma^{\text{pool}}$ is the same as in the pool market.

\begin{flalign}
	\min_{\Gamma^{\text{pool}}} \sum_{t \in \mathcal{T}} \sum_{n \in \mathcal{N}} \lambda^{\text{pool}}_{n,t} ( p^{\text{pool}}_{n,t}- z^{\text{pool}}_{n, t}) \label{equ:mix_market_obj2}
\end{flalign}

\subsection{Equivalent optimization problem}


The carbon market operator's problem is still the same as the problem under the P2P and pool market design.

Up to now, we have been using the same methodology for planning models that integrate the P2P market and the pool market. Due to this consistency, we are able to further investigate the mixed market by formulating an equivalent optimization problem. The objective function (\ref{equ:mix_central_obj}) aggregates the objectives of all the agents and is found to be the same as (\ref{equ:p2p_cen_obj}), since all the trading within the pool adds up to zero. 

{\allowdisplaybreaks
\begin{subequations}
\begin{flalign}
	\min_{\Gamma}
	 \sum_{n \in \mathcal{N}} f_n +  g + \sum_{t \in \mathcal{T}} \sum_{n \in \mathcal{N}} \sum_{m \in \omega_n} I_{n,m} |p^{\text{bi}}_{n,m,t}| \label{equ:mix_central_obj}
\end{flalign}
\begin{flalign}
    \text{s.t.} & (\ref{equ:mix_pro_balance1}) - (\ref{equ:mix_pro_others}), \forall n \in \mathcal{N} \label{equ:mix_central_constraint1} \\
    & (\ref{equ:mix_tso_ptdf}) - (\ref{equ:mix_tso_limit}) \label{equ:mix_central_constraint2}\\
	 & p^{\text{bi}}_{n,m,t} = - p^{\text{bi}}_{m,n,t}, \forall n \in \mathcal{N}, \forall m \in \omega_n, \forall t\in T  \text{: } \lambda^{\text{bi}}_{n,m,t} \label{equ:mix_central_constraint3} \\
	 & p^{\text{bi}}_{n,m,t} = z^{\text{bi}}_{n,m,t},  \forall n \in \mathcal{N}, \forall m \in \omega_n,\forall t\in T   \text{: } \lambda^{\text{grid}}_{n,m,t}  \label{equ:mix_central_constraint4}\\ 
 	 & p^{\text{pool}}_{n,t} - z^{\text{pool}}_{n,t} = 0, \qquad \forall n \in \mathcal{N}, \forall t\in T \text{: } \lambda^{\text{pool}}_{n,t} \label{equ:mix_central_constraint5} \\
	 & \sum_{n \in \mathcal{N}} e^{CO_2}_n = \text{CAP}^{CO_2} \text{: } \lambda^{CO_2} \label{equ:mix_cen_co2}
\end{flalign}
\end{subequations}
}

Constraints (\ref{equ:mix_central_constraint1}) and (\ref{equ:mix_central_constraint2}) are the same as those of the prosumers and the TSO. (\ref{equ:mix_central_constraint3}) - (\ref{equ:mix_central_constraint5}) are the KKT conditions from the optimization problems of the bilateral market operator and the pool market operator, where the dual variables of the constraints are the respective prices.

\subsection{Distributed solving formulation using ADMM}

(\ref{equ:mix_admm_L}) is the Lagrangian of (\ref{equ:mix_central_obj}). Compared to the Lagrangians of problems for the P2P market and the pool market alone, four coupling constraints are relaxed due to the co-existence of the two electricity markets.

{\allowdisplaybreaks
\begin{subequations}
\begin{flalign}
    \mathcal{L}(\Gamma^{\text{ps}}_n, \Gamma^{\text{TSO}}) = \sum_{n \in \mathcal{N}} f_n + g \notag \\
    + \sum_{t \in \mathcal{T}} \sum_{n \in \mathcal{N}} \sum_{m \in \omega_n} I_{n,m} |p^{\text{bi}}_{n,m,t}| \notag \\
     +	\sum_{t \in \mathcal{T}} \sum_{n \in \mathcal{N}} \sum_{m \in \omega_n} \lambda^{\text{bi}}_{n,m,t}  (p^{\text{bi}}_{n,m,t} + p^{\text{bi}}_{m,n,t}) \notag \\
     +	\sum_{t \in \mathcal{T}} \sum_{n \in \mathcal{N}} \sum_{m \in \omega_n} \frac{Q^{\text{bi}}_{n,m,t}}{2} (p^{\text{bi}}_{n,m,t} + p^{\text{bi}}_{m,n,t})^2 \notag \\
     + \sum_{t \in \mathcal{T}} \sum_{n \in \mathcal{N}} \sum_{m \in \omega_n} \lambda^{\text{grid}}_{n,m,t} (p^{\text{bi}}_{n,m,t} - z^{\text{bi}}_{n,m,t}) \notag \notag \\
     + \sum_{t \in \mathcal{T}} \sum_{n \in \mathcal{N}} \sum_{m \in \omega_n} \frac{Q^{\text{grid}}_{n,m,t}}{2} (p^{\text{bi}}_{n,m,t} - z^{\text{bi}}_{n,m,t})^2  \notag \\
    +	\sum_{t \in \mathcal{T}} \sum_{n \in \mathcal{N}}  \lambda^{\text{pool}}_{n,t} (p^{\text{pool}}_{n,t} - z^{\text{pool}}_{n, t})  \notag \\
    +	\sum_{t \in \mathcal{T}} \sum_{n \in \mathcal{N}}  \frac{Q^{\text{pool}}_{n,t}}{2}  (p^{\text{pool}}_{n,t} - z^{\text{pool}}_{n, t})^2 \notag \\
    +	\lambda^{CO_2} (\sum_{n \in \mathcal{N}} e^{CO_2}_n - \text{CAP}^{CO_2}) \notag \\ 
    +	\frac{Q^{CO_2}}{2}  (\sum_{n \in \mathcal{N}} e^{CO_2}_n - \text{CAP}^{CO_2})^2 \label{equ:mix_admm_L}
\end{flalign}

The ADMM algorithm is given as follows. The prosumers and the TSO solve their local optimization problems, respectively. The prices are updated by the two energy market operators and the carbon market operator.

\begin{flalign}
    \Gamma^{\text{ps}(k+1)}_n = \text{argmin}_{\Gamma^{\text{ps}}_n} \mathcal{L}(\Gamma^{\text{ps}}_n, \Gamma^{\text{TSO}(k)}),  \forall n \in \mathcal{N} \\
    \text{s.t. } (\ref{equ:mix_pro_balance1}) - (\ref{equ:mix_pro_others}) \notag \\    
    \Gamma^{\text{TSO}(k+1)} = \text{argmin}_{\Gamma^{\text{TSO}}} \mathcal{L}(\Gamma^{\text{ps}(k)}_n, \Gamma^{\text{TSO}})\\
    \text{s.t. } (\ref{equ:mix_tso_ptdf}) - (\ref{equ:mix_tso_limit}) \notag \\
    \lambda^{\text{bi}(k+1)}_{n,m,t} = \lambda^{\text{bi}(k)}_{n,m,t} +  Q^{\text{bi}}_{n,m,t} (p^{\text{bi}(k+1)}_{n,m,t} + p^{\text{bi}(k+1)}_{m,n,t}),  \notag\\
    \forall t \in \mathcal{T}, \forall n \in \mathcal{N}, \forall m \in \omega_n  \notag \\
    \lambda^{\text{grid}(k+1)}_{n,m,t} = \lambda^{\text{grid}(k)}_{n,m,t} +  Q^{\text{grid}}_{n,m,t} (p^{\text{grid}(k+1)}_{n,m,t} - z^{\text{bi}(k+1)}_{n,m,t}), \notag \\
    \forall t \in \mathcal{T}, \forall n \in \mathcal{N}, \forall m \in \omega_n \notag \\
    \lambda^{\text{pool}(k+1)}_{n,t} = \lambda^{\text{pool}(k)}_{n,t} + Q^{\text{pool}}_{n,t}  (p^{\text{pool}(k+1)}_{n,t} 
    - z^{\text{pool}(k+1)}_{n, t}), \notag \\
    \forall l \in \mathcal{L},  \forall t \in \mathcal{T} \notag \\
    \lambda^{{CO_2}(k+1)} = \lambda^{{CO_2}(k)} + Q^{CO_2}  (\sum_{n \in \mathcal{N}} e^{{CO_2}(k+1)}_n - \text{CAP}^{CO_2}) 
\end{flalign}

\end{subequations}

}

\section{Conclusion} \label{sec:conclusion}

In this paper, a prosumer-centric framework for concurrent generation and transmission planning is presented. Here, we consider prosumers as a general notion for an aggregated production and demand.
In current power system planning models, the market environment is often overlooked. Therefore, we consider the planning within three different electricity market environments, i.e. the P2P market with product differentiation, the pool market and the mixed bilateral/pool market.  We integrate the market models in the planning models, which combines and contributes to the two mainstream power system models, i.e. planning and market. Among all, one of our contributions is to consider the P2P market in a planning model whereas previous studies only focus on the operational aspects.  

For each of the three markets, the following modelling flow is used. We first formulate the planning optimization problems for different agents, i.e. the prosumers, the TSO, the electricity market operator and the carbon market operator, separately. This way of problem formulation exposes the true costs that are incurred for the agents along the planning horizon which include both the planning costs and the trading-related costs in the electricity market. 
After having formulated the optimization problems for all the agents, an equivalent centralized optimization problem has been found and then a distributed optimization algorithm is presented. 

Besides the P2P market, the pool market and the mixed bilateral/pool market are also included in the framework. These inclusions are particularly critical for two reasons. On the one hand, they allow investigating the influences of these markets on the planning decisions. On the other hand, as P2P markets prepare to roll out, the electricity market design might change in the coming decades, while investment decisions are made now under current market designs. The framework helps to assess the uncertainties related to the changing market designs. In addition, a carbon market with a cap-and-trade system is integrated into this study. It involves the government as an extra agent for the ease of modeling the effects of the carbon targets on the planning decisions, which is well on the political agenda all over the world.

Part II of the paper will continue with the numerical results from three archetypal case studies.

\bibliographystyle{IEEEtran}
\bibliography{IEEEabrv, ref}

\begin{thebibliography}{10}
\providecommand{\url}[1]{#1}
\csname url@samestyle\endcsname
\providecommand{\newblock}{\relax}
\providecommand{\bibinfo}[2]{#2}
\providecommand{\BIBentrySTDinterwordspacing}{\spaceskip=0pt\relax}
\providecommand{\BIBentryALTinterwordstretchfactor}{4}
\providecommand{\BIBentryALTinterwordspacing}{\spaceskip=\fontdimen2\font plus
\BIBentryALTinterwordstretchfactor\fontdimen3\font minus
  \fontdimen4\font\relax}
\providecommand{\BIBforeignlanguage}[2]{{%
\expandafter\ifx\csname l@#1\endcsname\relax
\typeout{** WARNING: IEEEtran.bst: No hyphenation pattern has been}%
\typeout{** loaded for the language `#1'. Using the pattern for}%
\typeout{** the default language instead.}%
\else
\language=\csname l@#1\endcsname
\fi
#2}}
\providecommand{\BIBdecl}{\relax}
\BIBdecl

\bibitem{Botterud2005}
A.~Botterrud, M.~D. Ilic, and I.~Wangensteen, ``{Optimal investments in power
  generation under centralized and decentralized decision making},'' \emph{IEEE
  Transactions on Power Systems}, vol.~20, no.~1, pp. 254--263, 2005.

\bibitem{Wang2017}
Y.~Wang, S.~Wang, and L.~Wu, ``{Distributed optimization approaches for
  emerging power systems operation: A review},'' \emph{Electric Power Systems
  Research}, vol. 144, 2017.

\bibitem{Parag2016}
Y.~Parag and B.~K. Sovacool, ``{Electricity market design for the prosumer
  era},'' \emph{Nature Energy}, vol.~1, no.~4, p. 16032, 2016.

\bibitem{VanLeeuwen2020}
G.~van Leeuwen, T.~AlSkaif, M.~Gibescu, and W.~van Sark, ``{An integrated
  blockchain-based energy management platform with bilateral trading for
  microgrid communities},'' \emph{Applied Energy}, vol. 263, no. October 2019,
  p. 114613, 2020.

\bibitem{Wang2020c}
N.~Wang, R.~Verzijlbergh, P.~Heijnen, and P.~Herder, ``{Modeling the
  decentralized energy investment and operation in the prosumer era: a
  systematic review},'' \emph{2020 IEEE PES Innovative Smart Grid Technologies
  Europe (ISGT-Europe)}, pp. 1079--1083, oct 2020.

\bibitem{Yang2015}
\BIBentryALTinterwordspacing
Y.~Yang, S.~Zhang, and Y.~Xiao, ``{Optimal design of distributed energy
  resource systems coupled withenergy distribution networks},'' \emph{Energy},
  vol.~85, pp. 433--448, 2015. [Online]. Available:
  \url{http://dx.doi.org/10.1016/j.energy.2015.03.101}
\BIBentrySTDinterwordspacing

\bibitem{Hahnel2020}
U.~J. Hahnel, M.~Herberz, A.~Pena-Bello, D.~Parra, and T.~Brosch, ``{Becoming
  prosumer: Revealing trading preferences and decision-making strategies in
  peer-to-peer energy communities},'' \emph{Energy Policy}, vol. 137, 2020.

\bibitem{Sorin2019}
E.~Sorin, L.~Bobo, P.~Pinson, and S.~Member, ``{Consensus-Based Approach to
  Peer-to-Peer Electricity Markets with Product Differentiation},'' \emph{IEEE
  Transactions on Power Systems}, vol.~34, no.~2, pp. 994--1004, 2019.

\bibitem{Morstyn2019a}
T.~Morstyn and M.~D. McCulloch, ``{Multiclass Energy Management for
  Peer-to-Peer Energy Trading Driven by Prosumer Preferences},'' \emph{IEEE
  Transactions on Power Systems}, vol.~34, no.~5, pp. 4005--4014, 2019.

\bibitem{Baroche2019b}
T.~Baroche, F.~Moret, and P.~Pinson, ``{Prosumer markets: A unified
  formulation},'' \emph{2019 IEEE Milan PowerTech, PowerTech 2019}, pp. 1--6,
  2019.

\bibitem{Baroche2019a}
T.~Baroche, P.~Pinson, R.~L.~G. Latimier, and H.~{Ben Ahmed}, ``{Exogenous Cost
  Allocation in Peer-to-Peer Electricity Markets},'' \emph{IEEE Transactions on
  Power Systems}, vol.~34, no.~4, pp. 2553--2564, 2019.

\bibitem{Bower2000}
J.~Bower and D.~W. Bunn, ``{Model-based comparisons of pool and bilateral
  markets for electricity},'' \emph{Energy Journal}, vol.~21, no.~3, 2000.

\bibitem{Knezevic2011}
G.~Kne{\v{z}}evi{\'{c}}, S.~Nikolovski, and P.~Maric, ``{Electricity spot
  market simulation involving bilateral contracts hedging},'' in \emph{2011 8th
  International Conference on the European Energy Market, EEM 11}, 2011.

\bibitem{Lopes2012}
F.~Lopes, T.~Rodrigues, and J.~Sousa, ``{Negotiating bilateral contracts in a
  multi-agent electricity market: A case study},'' in \emph{Proceedings -
  International Workshop on Database and Expert Systems Applications, DEXA},
  2012.

\bibitem{Imran2020}
K.~Imran, J.~Zhang, A.~Pal, A.~Khattak, K.~Ullah, and S.~M. Baig, ``{Bilateral
  negotiations for electricity market by adaptive agent-tracking strategy},''
  \emph{Electric Power Systems Research}, vol. 186, 2020.

\bibitem{Bompard2008}
E.~Bompard and Y.~Ma, ``{Modeling bilateral electricity markets: A complex
  network approach},'' \emph{IEEE Transactions on Power Systems}, vol.~23,
  no.~4, 2008.

\bibitem{Doorman2008}
G.~L. Doorman and A.~Botterud, ``{Analysis of generation investment under
  different market designs},'' \emph{IEEE Transactions on Power Systems},
  vol.~23, no.~3, 2008.

\bibitem{Jenabi2013}
M.~Jenabi, S.~M.~T. {Fatemi Ghomi}, and Y.~Smeers, ``{Bi-level game approaches
  for coordination of generation and transmission expansion planning within a
  market environment},'' \emph{IEEE Transactions on Power Systems}, vol.~28,
  no.~3, pp. 2639--2650, 2013.

\bibitem{Taheri2017}
S.~S. Taheri, J.~Kazempour, and S.~Seyedshenava, ``{Transmission expansion in
  an oligopoly considering generation investment equilibrium},'' \emph{Energy
  Economics}, vol.~64, pp. 55--62, 2017.

\bibitem{Grimm2020}
V.~Grimm, B.~R{\"{u}}ckel, C.~S{\"{o}}lch, and G.~Z{\"{o}}ttl, ``{The impact of
  market design on transmission and generation investment in electricity
  markets},'' \emph{Energy Economics}, p. 104934, 2020.

\bibitem{Tohidi2017}
Y.~Tohidi, L.~Olmos, M.~Rivier, and M.~R. Hesamzadeh, ``{Coordination of
  Generation and Transmission Development Through Generation Transmission
  Charges - A Game Theoretical Approach},'' \emph{IEEE Transactions on Power
  Systems}, vol.~32, no.~2, pp. 1103--1114, 2017.

\bibitem{Moret2020}
F.~Moret, A.~Tosatto, T.~Baroche, and P.~Pinson, ``{Loss allocation in joint
  transmission and distribution peer-to-peer markets},'' \emph{IEEE
  Transactions on Power Systems}, vol.~36, no.~3, pp. 1833--1842, 2020.

\bibitem{Hobbs03complementarity-basedequilibrium}
B.~F. Hobbs and U.~Helman, ``Complementarity-based equilibrium modeling for
  electric power markets,'' 2003.

\bibitem{StevenA.Gabriel2013}
{Steven A. Gabriel}, {Antonio J. Conejo}, {J. David Fuller}, {Benjamin F.
  Hobbs}, {Carlos Ruiz}, S.~A. Gabriel, A.~J. Conejo, J.~D. Fuller, B.~F.
  Hobbs, and C.~Ruiz, \emph{{Complementarity Modeling in Energy
  Markets}}.\hskip 1em plus 0.5em minus 0.4em\relax Springer-Verlag New York,
  2013, vol. 180.

\bibitem{Ruiz2014}
C.~Ruiz, A.~J. Conejo, J.~D. Fuller, S.~A. Gabriel, and B.~F. Hobbs, ``{A
  tutorial review of complementarity models for decision-making in energy
  markets},'' \emph{EURO Journal on Decision Processes}, vol.~2, no. 1-2, pp.
  91--120, 2014.

\bibitem{He2012}
Y.~He, L.~Wang, and J.~Wang, ``{Cap-and-trade vs. carbon taxes: A quantitative
  comparison from a generation expansion planning perspective},''
  \emph{Computers and Industrial Engineering}, vol.~63, no.~3, pp. 708--716,
  2012.

\bibitem{Molzahn2017}
D.~K. Molzahn, F.~D{\"{o}}rfler, H.~Sandberg, S.~H. Low, S.~Chakrabarti,
  R.~Baldick, and J.~Lavaei, ``{A Survey of Distributed Optimization and
  Control Algorithms for Electric Power Systems},'' \emph{IEEE Transactions on
  Smart Grid}, vol.~8, no.~6, 2017.

\bibitem{Kargarian2018}
A.~Kargarian, J.~Mohammadi, J.~Guo, S.~Chakrabarti, M.~Barati, G.~Hug, S.~Kar,
  and R.~Baldick, ``{Toward Distributed/Decentralized DC Optimal Power Flow
  Implementation in Future Electric Power Systems},'' \emph{IEEE Transactions
  on Smart Grid}, vol.~9, no.~4, 2018.

\bibitem{Sagastizabal2012}
C.~Sagastiz{\'{a}}bal, ``{Divide to conquer: Decomposition methods for energy
  optimization},'' in \emph{Mathematical Programming}, vol. 134, no.~1, 2012.

\bibitem{Ruiz2015}
C.~Ruiz and A.~J. Conejo, ``{Robust transmission expansion planning},''
  \emph{European Journal of Operational Research}, vol. 242, no.~2, 2015.

\bibitem{Munoz2016}
F.~D. Munoz, B.~F. Hobbs, and J.~P. Watson, ``{New bounding and decomposition
  approaches for MILP investment problems: Multi-area transmission and
  generation planning under policy constraints},'' \emph{European Journal of
  Operational Research}, vol. 248, no.~3, 2016.

\bibitem{Moret2019}
F.~Moret and P.~Pinson, ``{Energy Collectives: A Community and Fairness Based
  Approach to Future Electricity Markets},'' \emph{IEEE Transactions on Power
  Systems}, vol.~34, no.~5, pp. 3994--4004, 2019.

\bibitem{LeCadre2020}
H.~{Le Cadre}, P.~Jacquot, C.~Wan, and C.~Alasseur, ``{Peer-to-peer electricity
  market analysis: From variational to Generalized Nash Equilibrium},''
  \emph{European Journal of Operational Research}, vol. 282, no.~2, pp.
  753--771, 2020.

\bibitem{Wang2020b}
N.~Wang, R.~A. Verzijlbergh, P.~W. Heijnen, and P.~M. Herder, ``{A spatially
  explicit planning approach for power systems with a high share of renewable
  energy sources},'' \emph{Applied Energy}, vol. 260, p. 114233, feb 2020.

\end{thebibliography}

\end{document}